\documentclass{amsart}
\usepackage{amssymb}
\newtheorem{thm}{Theorem}
\newtheorem{lem}{Lemma}
\newtheorem{prop}{Proposition}
\newtheorem{cor}{Corollary}

\newtheorem{rem}{Remark}
\newtheorem{df}{Definition}

\begin{document}

\bibliographystyle{plain}

\title[Milo\v s Arsenovi\'c and Romi F. Shamoyan]{Embedding theorems for harmonic multifunction spaces on $\mathbb R^{n+1}_+$}

\author[]{Romi F. Shamoyan}
\author[]{Milo\v s Arsenovi\' c$\dagger$}

\address{Department of mathematics, University of Belgrade, Studentski Trg 16, 11000 Belgrade, Serbia}
\email{\rm arsenovic@matf.bg.ac.rs}

\address{Bryansk University, Bryansk Russia}
\email{\rm rshamoyan@yahoo.com}

\thanks{$\dagger$ Supported by Ministry of Science, Serbia, project OI174017}

\date{}

\begin{abstract}
We introduce and study properties of certain new multi functional harmonic spaces in the upper half space. We prove several sharp embedding theorems for such multi functional spaces, these results are new even in the case of a single function.
\end{abstract}

\maketitle

\footnotetext[1]{Mathematics Subject Classification 2010 Primary 42B15, Secondary 46E35.  Key words
and Phrases: multi functional spaces, harmonic functions, embedding theorems, upper half space.}

\section{Introduction and auxiliary results}

The theory of harmonic function spaces in the single function case is well developed by various authors during last decades, see \cite{St}, \cite{StW}, \cite{Gr1}. The main goal of this paper is to define, for the first time in the literature, multi functional harmonic spaces and to establish some properties of these spaces. The proofs we found and provided below are short modifications of proofs in the case of a single function, but even in this special case our results are new. We believe these new interesting objects can serve as a base for further generalizations and investigations.

Analytic analogues of theorems on multi functional spaces we obtained below were proved in recent papers of the first author \cite{LS1}, \cite{LS2}. We intend to expand these results to more general harmonic function spaces based on several functions in higher dimension. Let us note that these topics are new even in the classical case of harmonic function spaces on the unit disk.

We set $\mathbb H = \{(x, t) : x \in \mathbb R^n, t > 0 \} \subset \mathbb R^{n+1}$. For $z = (x, t) \in \mathbb H$ we set $\overline z = (x, -t)$. We denote the points in $\mathbb H$ usually by $z = (x, t)$ or $w = (y, s)$. The Lebesgue measure is denoted by $dm(z) = dz = dx dt$ or $dm(w) = dw = dy ds$, the Lebesgue measure of a measurable set $E \subset \mathbb H$ is often denoted by $|E|$. We also use measures $dm_\lambda(z) = t^\lambda dxdt$, $\lambda \in \mathbb R$.

The space of all harmonic functions in a domain $\Omega$ is denoted by $h(\Omega)$.
Weighted harmonic Bergman spaces on $\mathbb H$ are defined, for $0<p<\infty$ and $\lambda > -1$, by
$$A(p,\lambda) = A(p,\lambda)(\mathbb H) = \left\{ f \in h(\mathbb H) : \| f \|_{A(p,\lambda)} =
\left( \int_{\mathbb H} |f(z)|^p dm_\lambda(z) \right)^{1/p} < \infty \right\}.$$
These spaces are complete metric spaces, for $1\leq p < \infty$ they are Banach spaces.



For $f \in h(\mathbb H)$ and $t > 0$ we set $M_p(f, t) = \| f(\cdot, t) \|_{L^p(\mathbb R^n)}$, $0 < p <\infty$ with the
usual convention in the case $p = \infty$. We use harmonic mixed norm spaces
\begin{equation}
B^{p,q}_\alpha = \left\{ f \in h(\mathbb H) : \| f \|_{B^{p,q}_\alpha}^q = \int_0^\infty M_p^q(f, t) t^{\alpha q -1} dt
< \infty \right\},
\end{equation}
where $0\leq p < \infty$, $0<q<\infty$ and harmonic Triebel-Lizorkin spaces
\begin{equation}
F^{p,q}_\alpha = \left\{ f \in h(\mathbb H) : \| f \|_{F^{p,q}_\alpha}^p = \int_{\mathbb R^n} \left( \int_0^\infty
|f(x, t)|^q t^{\alpha q - 1} dt \right)^{p/q} dx < \infty \right\},
\end{equation}
where again $0 < p \leq \infty$ and  $0<q \leq \infty$, the case $q = \infty$ is covered by the usual convention. These spaces are complete metric spaces and for $\min(p,q) \geq 1$ they are Banach spaces. For details on $A(p,\lambda)$ spaces and more general $B^{p,q}_\alpha$ spaces see \cite{DS}; Triebel-Lizorkin spaces were studied, in the case of analytic functions, by many authors, see for example \cite{OF1}.

By $X$-$L^p$ Carleson measure (or simply $X$ Carleson measure when $p$ is clear from the context) of a (quasi)-normed subspace $X$ of $h(\Omega)$ we understand those positive Borel measures $\mu$ on $\Omega$ such that
\begin{equation}\label{decar}
\left( \int_\Omega |f|^p d\mu \right)^{1/p} \leq C \|f\|_X, \qquad f \in X.
\end{equation}
Typical cases are $\Omega = \mathbb B = \{x \in \mathbb R^n : |x| < 1 \}$ and $\Omega = \mathbb H$. In the case $\Omega$ is the unit disc and $X$ is the analytic Hardy space $H^p$ we have a classical notion of Carleson measures on the unit disc. The cases of harmonic spaces can be found, for example, in \cite{CKY} and much earlier in \cite{Sh1}.

One of the goals of this paper is to find complete descriptions of Carleson measures for certain new harmonic function spaces in higher dimension, and also in the multi function case, in the latter case some restrictions on the parameters involved
usually appear. We note that recently several new papers appeared where embedding theorems for analytic spaces in the unit disk were obtained and where the classical expression $\int_B|f|^pd\mu$  was modified, generalized or changed to certain expressions of $G(f,\mu,p)$, see for example \cite{GLW}, \cite{W1}; see also an earlier paper \cite{Co1}. In these papers descriptions of measures for which $G(f,\mu,p) \leq C \|f\|_X$ were presented. Our sharp embedding theorems 2, 3, 4 and 8 we present below are of this character, and the spaces we deal with are defined using the above mentioned idea. On the other hand, our theorems
1, 5 and 6 give results modeled after (\ref{decar}), but in multi functional case.


We use common convention regarding constants: letter $C$ denotes a constant which can change its value from one occurrence
to the next one. Given two positive quantities $A$ and $B$, we write $A \asymp B$ if there are two constants $c, C > 0$ such that $cA \leq B \leq CA$.






Many of the results of this paper rely on the following three key auxiliary results. The first one is a concrete Whitney
decomposition into cubes of the upper half space, the second one is essentially subharmonic behavior of $|f|^p$ for harmonic
$f$ and for $0<p<\infty$ and the third one describes geometric properties of Whitney cubes.

\begin{lem}[\cite{St}]\label{LemmaA}
There exists a collection $\{ \Delta_k \}_{k=1}^\infty$ of closed cubes in $\mathbb H$ with sides parallel to coordinate axes such that

$1^o$. $\cup_{k=1}^\infty \Delta_k = \mathbb H$ and ${\rm diam} \Delta_k \asymp {\rm dist} (\Delta_k,
\partial \mathbb H)$.

$2^o$.  The interiors of the cubes $\Delta_k$ are pairwise disjoint.

$3^o$.  If $\Delta_k^\ast$ is a cube with the same center as $\Delta_k$, but enlarged 5/4 times, then the
collection $\{ \Delta_k^\ast \}_{k=1}^\infty$ forms a finitely overlapping covering of $R^{n+1}_+$, i.e. there is
a constant $C = C_n$ such that $\sum_k \chi_{\Delta_k^\ast} \leq C$.
\end{lem}

\begin{lem}[\cite{D1}]\label{LemmaB}
Let $\Delta_k$ and $\Delta_k^\ast$ be the cubes from the previous lemma and let $(\xi_k, \eta_k)$ be the center of $\Delta_k$. Then, for $0<p<\infty$ and $\alpha > 0$, we have
\begin{equation}
\eta_k^{\alpha p - 1} \max_{z \in \Delta_k} |f(z)|^p \leq \frac{C}{|\Delta_k^\ast|} \int_{\Delta_k^\ast} t^{\alpha p - 1}
|f(x,t)|^p dx dt, \qquad f \in h(\mathbb H), \quad k \geq 1.
\end{equation}
\end{lem}

\begin{lem}[\cite{St}]\label{LemmaC}
Let $\Delta_k$ and $\Delta_k^\ast$ are as in the previous lemma, let $\zeta_k = (\xi_k, \eta_k)$ be the center of the cube $\Delta_k$. Then we have:
\begin{equation}\label{mlam}
 m_\lambda(\Delta_k) \asymp \eta_k^{n+1+\lambda} \asymp m_\lambda (\Delta_k^\ast), \qquad \lambda \in
\mathbb R,
\end{equation}
\begin{equation}\label{dis}
|\overline w - z | \asymp | \overline \zeta_k - z|, \qquad w \in \Delta_k^\ast, \quad z \in \mathbb H,
\end{equation}
\begin{equation}\label{height}
t \asymp \eta_k, \qquad (x,t) \in \Delta_k^\ast.
\end{equation}
\end{lem}

Since the cubes from the above Lemmata appear quite often in this paper, we fix the following notation:
$\Delta_k$ and $\Delta_k^\ast$ are the cubes from Lemma \ref{LemmaA}, centered at $\zeta_k = (\xi_k, \eta_k)$, $k \geq 1$.

We also need the following integral estimate.

\begin{lem}[\cite{KY}]\label{omit}
If $\alpha > -1$ and $n+\alpha<2\gamma-1$, then
\begin{equation}
\int_{\mathbb H} \frac{t^\alpha dz}{|z - \overline w|^{2\gamma}} \leq C s^{\alpha + n + 1 - 2\gamma}, \qquad
w = (y, s) \in \mathbb H.
\end{equation}
\end{lem}

The following definition introduces certain multi functional spaces that appeared in the setting of analytic functions in the
unit ball in $\mathbb C^n$ in \cite{LS3}.

\begin{df}
Let $s>-1$ and $\vec{p} = (p_1, \ldots, p_m)$ where $0 < p_i < \infty$. We denote by $A(\vec{p}, s, m)$ the set of all $m$-tuples $(f_1, \ldots, f_m)$ of functions harmonic in $\mathbb H$ such that
\begin{equation}\label{mulnor}
(f_1, \ldots, f_m)_{A(\vec{p}, s, m)} = \int_{\mathbb H} \prod_{i=1}^m |f_i(z)|^{p_i} t^s
dm(z) < \infty.
\end{equation}
If $p_i = p$, $i = 1, \ldots, m$, we write simply $A(p, s, m)$.
\end{df}

For the first part of the following proposition see \cite{AS5}, for the second one see \cite{AS4}.
Note that the first part gives a simple estimate of the $A(\vec{p}, \lambda, m)$ norm, while the second one gives an
estimate of trace in $A(p,\lambda)$ norm.

\begin{prop}\label{Lemma4}
$1^o$. Let $0<p_i<\infty$, $-1 < s_i < \infty$, $f_i \in A(p_i,s_i)(\mathbb H)$ for $i = 1, \ldots, m$ and set $\lambda =
(m-1)(n+1) + \sum_{i=1}^m s_i$. Then
\begin{equation}
(f_1, \ldots, f_m)_{A(\vec{p}, \lambda, m)}
\leq C \prod_{i=1}^m \| f_i \|_{A^(p_i,s_i)}^{p_i}.
\end{equation}
$2^o$. Let $0<p<\infty$ and $s_1, \ldots, s_m > -1$. Set $\lambda = (m-1)(n+1) + \sum_{j=1}^m s_j$. Then there is a constant
$C>0$ such that for all $f \in  h(\mathbb H^m)$ we have
\begin{equation}\label{eqL4}
\int_{\mathbb H} |f(z,\ldots, z)|^p dm_\lambda (z) \leq C \int_{\mathbb H} \cdots \int_{\mathbb H}
|f(z_1, \ldots, z_m)|^p dm_{s_1}(z_1) \ldots dm_{s_m}(z_m).
\end{equation}
\end{prop}



\section{Embedding theorems for multi functional spaces of harmonic functions}

In this section we give several sharp embedding theorems for multi functional spaces of harmonic functions, necessary and
sufficient conditions turn out to be Carleson-type conditions.

The following theorem is an analogue for harmonic functions spaces of Theorem 3 from \cite{LS1}.

\begin{thm}\label{Themb1}
Let $\mu$ a positive Boreal measure on $\mathbb H$. Assume $0<p_i, q_i < \infty$, $i = 1, \ldots m$, satisfy
$\sum_{i=1}^m \frac{1}{q_i} = 1$ and let $\alpha > -1$. Then the following two conditions on the measure $\mu$ are
equivalent.

$1^o$. If $f_i$, $i = 1, \ldots, m$ are functions harmonic in $\mathbb H$, then we have
\begin{equation}\label{Themeq1}
\int_{\mathbb H} \prod_{i=1}^m |f_i(z)|^{p_i} d\mu(z) \leq C \prod_{i=1}^m \left[ \sum_{k=1}^\infty \left(
\int_{\Delta_k^\ast} |f_i(z)|^{p_i} t^\alpha dz \right)^{q_i} \right]^{1/q_i}.
\end{equation}

$2^o$. The measure $\mu$ satisfies a Carleson type condition:
\begin{equation}\label{cacon1}
\mu(\Delta_k) \leq C |\Delta_k|^{m(1 + \frac{\alpha}{n+1})}, \qquad k \geq 1.
\end{equation}
\end{thm}

{\it Proof.} Let us assume $\mu$ satisfies condition (\ref{cacon1}). Then we have, using Lemma \ref{LemmaA} and Lemma \ref{LemmaB}:
\begin{align*}
\phantom{M} & \int_{\mathbb H} \prod_{i=1}^m |f_i(z)|^{p_i} d\mu(z) = \sum_{k=1}^\infty \int_{\Delta_k} \prod_{i=1}^m |f_i(z)|^{p_i} d\mu(z) \leq \sum_{k=1}^\infty \mu(\Delta_k) \prod_{i=1}^m \sup_{\Delta_k} |f_i|^{p_i}\\
& \leq C\sum_{k=1}^\infty \mu(\Delta_k) \eta_k^{-m(n+1+\alpha)}
\prod_{i=1}^m \int_{\Delta_k^\ast} |f_i(w)|^{p_i} s^\alpha dw\\
& \leq C\sum_{k=1}^\infty \prod_{i=1}^m \int_{\Delta_k^\ast} |f_i(w)|^{p_i} s^\alpha dw.
\end{align*}
Now one arrives at estimate (\ref{Themeq1}) by applying Holder's inequality for sums with exponents $q_1, \ldots, q_m$.

Now we assume (\ref{Themeq1}) holds. We fix $k \in \mathbb N$ and choose $f_i(z) = f(z) = |z - \overline \zeta_k|^{1-n}$. Clearly
$$\prod_{i=1}^m |f_i(z)|^{p_i} = |z - \overline\zeta_k|^{-(n-1)\sum_{i=1}^m p_i}.$$
Therefore (\ref{Themeq1}), combined with Lemma \ref{LemmaC}, gives
\begin{align*}
\frac{\mu(\Delta_k)}{\eta_k^{(n-1)\sum_{i=1}^m p_i}} &\leq C \int_{\Delta_k} |z - \overline\zeta_k|^{-(n-1)\sum_{i=1}^m p_i}
d\mu(z) \leq C \int_{\mathbb H} \prod_{i=1}^m |f_i(z)|^{p_i} d\mu(z) \\
& \leq C \prod_{i=1}^m \int_{\Delta_k^\ast} \frac{t^\alpha dz}{|z -
\overline \zeta_k|^{p_i(n-1)}} \leq C \prod_{i=1}^m |\Delta_k| \eta_k^{\alpha - p_i(n-1)},
\end{align*}
and (\ref{cacon1}) easily follows. $\Box$



The spaces defined below were considered, in the case of analytic functions on the unit ball in $\mathbb C^n$, in \cite{LS2}.

\begin{df}
Let $0<p,q<\infty$ and let $\mu$ be a positive Borel measure on $\mathbb H$. The space $A(p,q,m,d\mu)$ is the space of all $(f_1, \ldots, f_m)$ where $f_i \in h(\mathbb H)$ for $i=1, \ldots, m$ such that
\begin{equation}\label{apq}
\| (f_1, \ldots, f_m) \|_{A(p,q,m,d\mu)}^q = \sum_{k=1}^\infty \left( \int_{\Delta_k}
\prod_{i=1}^m |f_i(z)|^p d\mu(z) \right)^{q/p} < \infty.
\end{equation}
If $d\mu = dm_s$, then we write $A(p,q,m,s)$ for the corresponding space, if $m=1$ we write $A(p,q,d\mu)$ or $A(p,q,s)$ if
$d\mu = dm_s$.
\end{df}


\begin{lem}\label{isum}
Let $0<s<\infty$ and $\beta > -1$. Then we have
\begin{equation}\label{isumeq}
\int_{\mathbb H} |f(z)|^s t^\beta dx dt \asymp \sum_{k=1}^\infty \eta_k^{n+1+\beta} \sup_{\Delta_k} |f|^s, \qquad
f \in h(\mathbb H).
\end{equation}
\end{lem}

We omit an easy proof based on Lemma \ref{LemmaA} and Lemma \ref{LemmaC}.

Theorem 3.1 from \cite{LS2} served as a model for the following theorem.

\begin{thm}\label{emapq}
Let $0<p,q<\infty$, $0<s\leq q$, $\beta_i > -1$ for $i = 1, \ldots, m$ and let $\mu$ be a positive Borel measure on $\mathbb H$. Then the following conditions are equivalent.

$1^o$. If $f_i \in A(s,\beta_i)$, $i = 1, \ldots, m$, then
\begin{equation}\label{emapqeq}
\| (f_1, \ldots, f_m) \|_{A(p,q,m,d\mu)} \leq C \prod_{i=1}^m \| f \|_{A(s, \beta_i)}.
\end{equation}

$2^o$. The measure $\mu$ satisfies a Carleson type condition:
$$\mu(\Delta_k) \leq C \eta_k^{\sum_{i=1}^m \frac{p(n+1+\beta_i)}{s}}, \qquad k \geq 1.$$
\end{thm}

{\it Proof.} Assume $2^o$ holds and choose $f_i \in A(s, \beta_i)$, $i = 1, \ldots, m$. Since $0<s/q\leq 1$ we have, using Lemma \ref{isum}
\begin{align*}
\| (f_1, \ldots, f_m) \|_{A(p,q,m,d\mu)}^s & = \left( \sum_{k=1}^\infty \left( \int_{\Delta_k} \prod_{i=1}^m
|f_i(z)|^p d\mu(z) \right)^{q/p} \right)^{s/q}\\
& \leq \left( \sum_{k=1}^\infty \mu(\Delta_k)^{q/p} \max_{\Delta_k} \prod_{i=1}^m |f_i|^q \right)^{s/q}\\
& \leq C \sum_{k=1}^\infty \prod_{i=1}^m \eta_k^{n+1+\beta_i} \max_{\Delta_k} |f_i|^s \\
& \leq C \prod_{i=1}^m \left( \sum_{k_i = 1}^\infty \eta_k^{n+1+\beta_i} \max_{\Delta_k} |f_i|^s \right) \leq C \prod_{i=1}^m \| f_i \|_{A(s,\beta_i)}^s,
\end{align*}
and we proved implication $2^o \Rightarrow 1^o$.

Conversely, assume $1^o$ holds and choose $l \in \mathbb N_0$ such that $s(n+l-1) > n + \beta + 1$. We use, as test functions,
functions $f_{\zeta_k,l}$, $k \in \mathbb N$ where
\begin{equation}\label{testf}
f_{w,l}(z) = \frac{\partial^l}{\partial t^l} \frac{1}{|z-\overline w|^{n-1}}, \qquad z \in \mathbb H,
\end{equation}
see \cite{AS4} for norm and pointwise estimates related to these functions. In particular we have
\begin{equation}
\| f_{\zeta_k,l} \|_{A(s,\beta_i)} \leq C \eta_k^{-(n-1+l) + \frac{n+\beta_i+1}{s}}, \qquad k \geq 1, \quad 1\leq i \leq m
\end{equation}
and
\begin{align}
\| (f_{\zeta_k,l}, \ldots, f_{\zeta_k,l}) \|_{A(p,q,m,d\mu)} & \geq \left( \int_{\Delta_k} |f_{\zeta_k,l}(z)|^{mp} d\mu(z) \right)^{1/p} \notag \\
& \geq C \mu(\Delta_k)^{1/p} \eta_k^{-m(n+l-1)}.
\end{align}
The last two estimates combined with (\ref{emapqeq}) give $2^o$. $\Box$

The corresponding result for mixed norm spaces is the following theorem.

\begin{thm}\label{emapqf}
Let $0<p,q<\infty$, $0<s\leq q$, $\beta_i > -1$ for $i = 1, \ldots, m$, $0 < t_i \leq s$ for $i = 1, \ldots, m$ and let $\mu$ be a positive Borel measure on $\mathbb H$. Then the following conditions are equivalent:

$1^o$. If $f_i \in B^{s, t_i}_{(\beta_i + 1)/s}$, $i = 1, \ldots, m$, then
\begin{equation}\label{emapqeq}
\| (f_1, \ldots, f_m) \|_{A(p,q,m,d\mu)} \leq C \prod_{i=1}^m \| f \|_{B^{s,t_i}_{\frac{\beta_i + 1}{s}}}.
\end{equation}

$2^o$. If $f_i \in F^{s, t_i}_{(\beta_i + 1)/s}$, $i = 1, \ldots, m$, then
\begin{equation}\label{emapqeq}
\| (f_1, \ldots, f_m) \|_{A(p,q,m,d\mu)} \leq C \prod_{i=1}^m \| f \|_{F^{s,t_i}_{\frac{\beta_i + 1}{s}}}.
\end{equation}

$3^o$. The measure $\mu$ satisfies a Carleson type condition:
$$\mu(\Delta_k) \leq C \eta_k^{\sum_{i=1}^m \frac{p(n+1+\beta_i)}{s}}, \qquad k \geq 1.$$
\end{thm}

{\it Proof.} Since $A(s,\beta) = B^{s,s}_{(\beta + 1)/s}$, the previous theorem, combined with embeddings
$B^{s,t}_{(\beta + 1)/s} \hookrightarrow B^{s,s}_{(\beta + 1)/s}$ and $F^{s,t}_{(\beta + 1)/s} \hookrightarrow
B^{s,s}_{(\beta + 1)/s}$, valid for $0 < t \leq s$, gives implications $3^o \Rightarrow 1^o$ and $3^o \Rightarrow 2^o$.

Now we assume $2^o$ holds. Let us choose $l \in \mathbb N_0$ such that $s(n+l-1) > n + \beta + 1$. Set, for $w = (y, \sigma)
\in \mathbb H$, $f_i(z) = f(z) = f_{w,l}(z)$, $i = 1, \ldots, m$, see (\ref{testf}). We have, see \cite{AS4}:
\begin{equation}\label{as4}
\| f_{w,l} \|_{F^{s,t}_{\frac{\beta + 1}{s}}} = C \sigma^{\frac{n}{s} - (n-1+l - \frac{\beta+1}{s})}.
\end{equation}
Let us fix a cube $\Delta_k$. Using pointwise estimates from below for functions $f_{w,l}$, see \cite{AS4}, we obtain:
\begin{align}
\mu(\Delta_k)^{1/p} \eta_k^{-m(n-1+l)} & \leq C \left(\int_{\Delta_k} |f_{w,l}(z)|^{mp} d\mu(z) \right)^{1/p} \notag \\
& \leq C \| (f_1, \ldots, f_m) \|_{A(p,q,m,d\mu)},   \label{fpq} \qquad k \geq 1.
\end{align}
Now combining (\ref{fpq}), (\ref{as4}) with $w = \zeta_k$ and (\ref{emapqeq}) gives $\mu(\Delta_k) \leq C
\eta_k^{\sum_{i=1}^m \frac{p(n+1+\beta_i)}{s}}$ and $3^o$ follows.

The implication $1^o \Rightarrow 3^o$ can be proved using the same test functions $f_{w,l}$ as above, relevant estimates
of $B^{p,q}_\alpha$ norm of these functions can be found in \cite{AS4} and we leave details to the reader. $\Box$

For $w = (y, s) \in \mathbb H$ we set $Q_w$ to be the cube, with sides parallel to the coordinate axis, centered at $w$ with
side length equal to $s$.

Our next multi functional embedding theorem is motivated by Theorem 3.6 from \cite{LS2}, which deals with the case of
a single analytic function on the unit ball in $\mathbb C^n$.

\begin{thm}\label{trisest}
Let $0<p,q < \infty$ and $0<\sigma_i \leq q$, $-1 < \alpha_i < \infty$ for $i = 1, \ldots, m$. Let $\mu$ be a positive Borel
measure on $\mathbb H$. Then the following two conditions are equivalent:

$1^o$. For any $m$-tuple $(f_1, \ldots, f_m)$ of harmonic functions on $\mathbb H$ we have
\begin{equation}\label{trsteq}
\| (f_1, \ldots, f_m) \|_{A(p,q,m,d\mu)}^q \leq C \prod_{i=1}^m \int_{\mathbb H} \left( \int_{Q_w}
|f_i(z)|^{\sigma_i} t^{\alpha_i}
dz \right)^{q/\sigma_i} dw.
\end{equation}

$2^o$. The measure $\mu$ satisfies the following Carleson-type condition:
\begin{equation}\label{cacoq}
\mu(\Delta_k) \leq C \eta_k^{m(n+1)\frac{p}{q} + \sum_{i=1}^m \frac{p(n+1+\alpha_i)}{\sigma_i}}, \qquad k \geq 1.
\end{equation}
\end{thm}

{\it Proof.} Let us prove $2^o \Rightarrow 1^o$ assuming $m=1$. Let $f \in h(\mathbb H)$, then we have, using Lemma
\ref{LemmaA}, Lemma \ref{LemmaB} and Lemma \ref{LemmaC}:

\begin{align*}
\| f \|_{A(p,q,d\mu)}^q & = \sum_{k=1}^\infty \left( \int_{\Delta_k} |f(z)|^p d\mu(z) \right)^{q/p} \leq \sum_{k=1}^\infty
\left( \max_{\Delta_k} |f|^p \mu(\Delta_k) \right)^{q/p}\\
& = \sum_{k=1}^\infty \max_{\Delta_k} |f|^q \eta_k^{q(\frac{n+1+\alpha}{\sigma} + \frac{n+1}{q})}\\
& \leq C \sum_{k=1}^\infty \left( \int_{\Delta_k^\ast} |f(z)|^\sigma t^{-n-1} dz \right)^{q/\sigma}
\eta_k^{q(\frac{n+1+\alpha}{\sigma} + \frac{n+1}{q})}\\
& \leq C \sum_{k=1}^\infty \left( \int_{\Delta_k^\ast} |f|^\sigma dz \right)^{q/\sigma} \eta_k^{n+1 + q\alpha/\sigma}.
\end{align*}
We continue this estimation, using Lemma \ref{LemmaB} and Lemma \ref{LemmaC}, and obtain
\begin{align}
\| f \|_{A(p,q,d\mu)}^q & \leq C \sum_{k=1}^\infty \left( \int_{\Delta_k^\ast} \left( \int_{Q_z} |f(w)|^\sigma s^\alpha dw
\right) t^{-n-1-\alpha} dz \right)^{q/\sigma} \eta_k^{n+1+ q\alpha/\sigma} \notag \\
& \leq C \sum_{k=1}^\infty \left( \int_{\Delta_k^\ast} \left( \int_{Q_z} |f(w)|^\sigma s^\alpha dw \right) t^{-n-1} dz
\right)^{q/\sigma} \eta_k^{n+1} \notag \\
& \leq C \sum_{k=1}^\infty \int_{\Delta_k^\ast} \left( \int_{Q_z} |f(w)|^\sigma s^\alpha dw \right)^{q/\sigma} t^{-n-1} dz \;
\eta_k^{n+1} \label{holder} \\
& \leq C \sum_{k=1}^\infty \int_{\Delta_k^\ast} \left( \int_{Q_z} |f(w)|^\sigma s^\alpha dw \right)^{q/\sigma} dz \notag \\
& \leq C \int_{\mathbb H} \left( \int_{Q_z} |f(w)|^\sigma s^\alpha dw \right)^{q/\sigma} dz \notag,
\end{align}
arriving at (\ref{trsteq}) for $m=1$. Note that at (\ref{holder}) we used Holder's inequality and at the last step we used finite overlapping property of the cubes $\Delta_k^\ast$.

Now we assume $1^o$, and again we restrict ourselves to the case $m=1$. Let us choose $l \in \mathbb N_0$ such that
$q(n-1+l) > q(n+1+\alpha)/\sigma + n + 1$. Let us fix a cube $\Delta_k$ centered at $\zeta_k = (\xi_k, \eta_k)$. We use a
construction from \cite{AS4}, where interested reader can find more details. Namely, there are constants $c> 0$ and $\delta > 0$ such that for all $w \in \mathbb H$ we have
\begin{equation}
|T_w| = c|Q_w|, \qquad T_w = \{ z \in Q_w : |f_{w,l}(z)| > c |z-\overline w|^{-(n-1+l)} \}.
\end{equation}
Next, using a compactness argument, one shows that for $w=(y,s) \in \mathbb H$ there are points $w_j = (y_j, s_j) \in
\mathbb H$, $1\leq j \leq N$, such that $s_j \asymp s$ and $Q_w = \cup_{j=1}^N T_{w_j} \cap Q_w$. Here $N$ depends only on
$n$ and $l$. We apply this argument to $w = \zeta_k$, it easily follows that there is $\theta_k = (u_k, v_k) \in \mathbb H$
such that $\mu(T_{\theta_k} \cap \Delta_k) \geq N^{-1} \mu(\Delta_k)$ and $v_k \asymp \eta_k$. Now we choose $f = f_{\theta_k, l}$ as a test function. We have
\begin{align}
\| f \|_{A(p,q,d\mu)}^q & \geq \left( \int_{\Delta_k} |f(z)|^p d\mu(z) \right)^{q/p} \geq
\left( \int_{\Delta_k\cap T_{\theta_k}} |f(z)|^p d\mu(z) \right)^{q/p} \notag \\
& \geq C \mu(\Delta_k \cap T_{\theta_k})^{q/p} v_k^{-q(n-1+l)} \notag \\
& \geq C \mu(\Delta_k) \eta_k^{-q(n-1+l)}. \label{esbl}
\end{align}
On the other hand, using Lemma \ref{LemmaC} and Lemma \ref{omit}, we have
\begin{align*}
\int_{\mathbb H} \left( \int_{Q_w} |f(z)|^\sigma dm_\alpha (z) \right)^{q/\sigma} dw & \leq C \int_{\mathbb H}
\left( \frac{s^{n+1+\alpha}}{|w-\overline \theta_k |^{\sigma(n-1+l)}} \right)^{q/\sigma} dw\\
& \leq C v_k^{\frac{q}{\sigma}(n+1+\alpha) - q(n-1+l) + n + 1}\\
& \leq C \eta_k^{\frac{q}{\sigma}(n+1+\alpha) - q(n-1+l) + n + 1}.
\end{align*}
This estimate, combined with (\ref{esbl}) and (\ref{trsteq}) (with $m=1$), gives desired estimate (\ref{cacoq}) for
$m=1$.

The general case $m>1$ can be proved along the same lines, multiple sums as in the proof of Theorem \ref{emapq} appear. Since no new ideas are involved we can leave details to the interested reader. $\Box$

As a preparation for our next result we state and prove the following lemma.

\begin{lem}\label{prepl}
Let $\alpha > -1$, $0< \tau < \infty$. Then there is a constant $C$ such that for any measurable function $u(z) \geq
0$ on $\mathbb H$ we have
\begin{equation}\label{eqprepl}
\eta_k^{n+1} \left( \int_{\Delta_k^\ast} u(z) dm_{\alpha}(z) \right)^\tau \leq C \int_{\Delta_k^\ast} \left(
\int_{Q_w} u(z) dm_\alpha (z) \right)^\tau dw, \quad k \geq 1.
\end{equation}
\end{lem}

{\it Proof.} Let us denote the left hand side by $L_\alpha$ and the right hand side, without constant $C$, by $D_\alpha$.
Since $L_\alpha \asymp \eta_k^\alpha L_0$ and $D_\alpha \asymp \eta_k^\alpha D_0$ it suffices to consider the special case
$\alpha = 0$. Let $q_w$ denote the cube centered at $w = (y, s) \in \mathbb H$, with side length equal $4s/5$. Then we have
$\Delta_k^\ast = \cup_{i=1}^N q_{w_i}$ for suitable $w_1, \ldots, w_N \in \Delta_k^\ast$ where $N$ depends only on $n$.
Therefore (\ref{eqprepl}) reduces to the following estimate:
\begin{equation}\label{redpr}
\eta_k^{n+1} \left( \int_{q_w} u(z) dm(z) \right)^\tau \leq C \int_{\Delta_k^\ast} \left(
\int_{Q_{\tilde w}} u(z) dm(z) \right)^\tau d\tilde w, \qquad w \in \Delta_k^\ast.
\end{equation}
Now we fix $w \in \Delta_k^\ast$. Clearly, $\int_{q_w} u(z) dm(z) \leq \int_{Q_{\tilde w}} u(z) dm(z)$ whenever $q_w \subset Q_{\tilde w}$, and the estimate (\ref{redpr}) follows from the simple observation, based on the sizes of $q_w$ and $Q_w$, that
\begin{equation*}
|E_k| \geq c \eta_k^{n+1}, \qquad E_k = \{ \tilde w \in \Delta_k^\ast : q_w \subset Q_{\tilde w} \}. \quad \Box
\end{equation*}

In the following theorem we allow for different exponents $p_i$.

\begin{thm}\label{trised}
Let $0<p_i, \sigma_i<\infty$, $-1 < \alpha_i < \infty$ for $i = 1, \ldots, m$. Let $\mu$ be a positive Borel measure on $\mathbb H$. Then the following two conditions are equivalent:

$1^o$. For any $m$-tuple $(f_1, \ldots, f_m)$ of harmonic functions on $\mathbb H$ we have
\begin{equation}\label{trsteqn}
\int_{\mathbb H} \prod_{i=1}^m|f_i(z)|^{p_i} d \mu(z) \leq C \prod_{i=1}^m \int_{\mathbb H}
\left( \int_{Q_w} |f_i(z)|^{\sigma_i} dm_{\alpha_i}(z) \right)^{p_i/\sigma_i} dw .
\end{equation}

$2^o$. The measure $\mu$ satisfies the following Carleson-type condition:
\begin{equation}\label{cacoqn}
\mu(\Delta_k) \leq C \eta_k^{m(n+1) + \sum_{i=1}^m \frac{p_i(n+1+\alpha_i)}{\sigma_i}},
\qquad k \geq 1.
\end{equation}
\end{thm}

{\it Proof.} Let us assume (\ref{trsteqn}) holds. We choose $l \in \mathbb N$ such that
$$p_i (n-1+l) > n+1+ \frac{p_i}{\sigma_i} (n + 1 + \alpha_i) \qquad  i = 1, \ldots, m.$$
We fix a cube $\Delta_k$ and take as test functions $f_i = f = f_{\theta_k, l}$, $i = 1, \ldots, m$, where $\theta_k =
(u_k, v_k)$ is chosen as in Theorem \ref{trisest}. Then we have, using Lemma \ref{omit}:
\begin{align*}
\eta_k^{-(n-1+l) \sum_{i=1}^m p_i} \mu (\Delta_k)& \leq C \int_{T_{\theta_k} \cap \Delta_k}
\prod_{i=1}^m|f_i(z)|^{p_i} d \mu(z)
\leq C \int_{\mathbb H} \prod_{i=1}^m|f_i(z)|^{p_i} d \mu(z)\\
& \leq C \prod_{i=1}^m \int_{\mathbb H} \left( \int_{Q_w} |f_i(z)|^{\sigma_i} dm_{\alpha_i}(z) \right)^{p_i/\sigma_i} dw\\
& \leq C \prod_{i=1}^m \int_{\mathbb H}\left( \frac{s^{n+1+\alpha_i}}{|\overline \theta_k - w|^{\sigma_i (n-1+l)}}
\right)^{p_i/\sigma_i} dw \\
& \leq  C v_k^{m(n+1) + \sum_{i=1}^m \frac{p_i(n+1+\alpha_i)}{\sigma_i} - (n-1+l)\sum_{i=1}^m p_i}\\
& \leq  C \eta_k^{m(n+1) + \sum_{i=1}^m \frac{p_i(n+1+\alpha_i)}{\sigma_i} - (n-1+l)\sum_{i=1}^m p_i}
\end{align*}
which gives (\ref{cacoqn}). Conversely, we assume now (\ref{cacoqn}). Using finite overlapping property of $\Delta_k^\ast$
we obtain
\begin{equation*}
\prod_{i=1}^m \int_{\mathbb H}
\left( \int_{Q_w} |f_i(z)|^{\sigma_i} dm_{\alpha_i}(z) \right)^{\frac{p_i}{\sigma_i}} dw \geq C \sum_{k=1}^\infty
\prod_{i=1}^m \int_{\Delta_k^\ast} \left( \int_{Q_w} |f_i(z)|^{\sigma_i} dm_{\alpha_i}(z) \right)^{\frac{p_i}{\sigma_i}} dw
\end{equation*}
and we have also an obvious estimate
\begin{align*}
\int_{\Delta_k} \prod_{i=1}^m|f_i(z)|^{p_i} d \mu(z) & \leq \mu(\Delta_k) \prod_{i=1}^m \max_{z \in \Delta_k} |f_i(z)|^{p_i}\\
& \leq C \eta_k^{m(n+1) + \sum_{i=1}^m \frac{p_i(n+1+\alpha_i)}{\sigma_i}}\prod_{i=1}^m \max_{z \in \Delta_k} |f_i(z)|^{p_i}.
\end{align*}
Therefore it suffices to prove, for each $k \geq 1$ and $1 \leq i \leq m$:
\begin{equation*}
\eta_k^{n+1+\frac{p_i(n+1+\alpha_i)}{\sigma_i}} \max_{z \in \Delta_k} |f_i(z)|^{p_i} \leq C \int_{\Delta_k^\ast}
\left( \int_{Q_w} |f_i(z)|^{\sigma_i} dm_{\alpha_i}(z) \right)^{p_i/\sigma_i} dw.
\end{equation*}
However, for $k \geq 1$ and $1 \leq i \leq m$, using Lemma \ref{LemmaB} and Lemma \ref{prepl} we obtain:
\begin{align*}
\eta_k^{n+1+\frac{p_i(n+1+\alpha_i)}{\sigma_i}} \max_{z \in \Delta_k} |f_i(z)|^{p_i}& \leq C \eta_k^{n+1} \left(
\int_{\Delta_k^\ast} |f_i(z)|^{\sigma_i} dm_{\alpha_i}(z) \right)^{\frac{p_i}{\sigma_i}}\\
& \leq C \int_{\Delta_k^\ast} \left( \int_{Q_w} |f(z)|^{\sigma_i} dm_{\alpha_i}(z) \right)^{\frac{p_i}{\sigma_i}} dw
\end{align*}
and the proof is completed. $\Box$


The following theorem is motivated by Theorem A from \cite{LS3}. It gives a very general and flexible method of obtaining
multi functional results starting from a uniform estimate for a single function.

\begin{thm}\label{thmA}
Let $\mu$ be a positive Borel measure on an open proper subset $G$ of $\mathbb R^n$ and let $(X_i, \|\cdot \|_{X_i})$ be a
(quasi)-normed space of functions harmonic in $G$, $1\leq i \leq k$. Let $d(x) = {\rm dist}(x, \partial G)$, $x \in G$ and
$0<q_i<\infty$, $-1 < \beta_i < \infty$ for $i = 1, \ldots, k$. Assume
\begin{equation}\label{unest}
\sup_{x\in G} |f_i(x)|^{q_i} d(x)^{\beta_i} \leq C \| f_i \|_{X_i}^{q_i}, \qquad f_i \in X_i, \quad i = 2, \ldots, k
\end{equation}
and
\begin{equation}
\int_G |f_1(x)|^{q_1} d\mu(x) \leq C \| f_1 \|_{X_1}^{q_1}, \qquad f_1 \in X_1.
\end{equation}
Then
\begin{equation}
\int_G \prod_{i=1}^k |f_i(x)|^{q_i} d(x)^{\beta_2 + \cdots + \beta_k} d\mu(x) \leq C \prod_{i=1}^k \| f_i \|_{X_i}^{q_i},
\qquad f_i \in X_i.
\end{equation}
\end{thm}

We refer the reader for the proof by induction to \cite{LS3}, where arguments readily extend to the present situation.
We also note that the above theorem can be extended to weights more general than $d(x)$, however it is precisely the
case of $d(x)$ where uniform estimates of type (\ref{unest}) are available. For example, if $G = \mathbb H$, then
$d(z) = t$ and we have an estimate
\begin{equation*}
|f(z)|^p d(z)^{n+1+\alpha} \leq C \| f \|_{A(p,\alpha)}^p, \qquad f \in A(p,\alpha), \quad z \in \mathbb H.
\end{equation*}
Therefore, taking $d\mu(z) = t^{s_1} dxdt$ we obtain the following corollary.

\begin{cor}
If $0<p_i<\infty$, $-1 < s_i < \infty$, $f_i \in A^(p_i,s_i)(\mathbb H)$ for $i = 1, \ldots, m$ and $\alpha > -1$, then
\begin{equation}
\int_{\mathbb H} \prod_{i=1}^m |f_i(z)|^{p_i} t^{(n+1)(m-1) + \sum_{i=1}^m s_i} dz \leq C \prod_{i=1}^m
\| f_i \|_{A(p_i,s_i)}^{p_i}
\end{equation}
\end{cor}

Using embedding $B^{p,q}_\alpha \hookrightarrow B^{p,p}_\alpha$, $0<q\leq p < \infty$, see \cite{DS}, we deduce another corollary.

\begin{cor}
If $0<q_i\leq p_i < \infty$, $-1 < s_i < \infty$ and $f_i \in B^{p_i, q_i}_{s_i}$ for $1 \leq i \leq m$, then
\begin{equation}
\int_{\mathbb H} \prod_{i=1}^m |f_i(z)|^{p_i} t^{n(m-1) + \sum_{i=1}^m s_ip_i -1} dz \leq C \prod_{i=1}^m
\| f_i \|_{B^{p_i, q_i}_{s_i}}^{p_i}.
\end{equation}
\end{cor}

Let $\mathcal D^\beta f$ denote the Riesz potential of $f \in h(\mathbb H)$ with respect to $t$ variable. Using estimate
\begin{equation}
\int_{\mathbb H} |f(z)|^p t^{\alpha p -1} dz \leq C \int_{\mathbb H} |\mathcal D^\beta f(z)|^p t^{\alpha p -1 + \beta p} dz
= N_{\alpha, \beta, p}(f), \qquad f \in h(\mathbb H)
\end{equation}
which is valid for $0<p<\infty$, $\alpha, \beta > 0$, see Theorem 7 from \cite{K1}, we obtain the following corollary.

\begin{cor}
Let $0<s_i, p_i < \infty$ and $f_i \in A(p_i,s_i)$ for $i=1, \ldots, m$. Then we have
\begin{equation}
\int_{\mathbb H} \prod_{i=1}^m |f_i(z)|^{p_i} t^{(n+1)(m-1) + \sum_{i=1}^m s_i} dz
\leq C \prod_{i=1}^m N_{s_i, \beta, p_i}(f_i).
\end{equation}
\end{cor}
Similar assertions can be obtained easily also in the unit ball $\mathbb B$.

Let $H^p(\mathbb B)$, $1<p<\infty$, stand for the harmonic Hardy space on the open unit ball.

\begin{prop}
Let $f_i \in H^{p_i}(\mathbb B)$ where $1<p_i<\infty$ for $i = 1, \ldots, m$. Then we have
\begin{equation}
\int_{\mathbb S} \prod_{j=1}^m |f_j(r \xi)|^{p_j} d\sigma(\xi)\;(1-r)^{(m-1)n} \leq C \prod_{j=1}^m
\| f_j \|_{H^{p_j}}^{p_j}.
\end{equation}
Similarly, if $f_i \in H^{p_i}(\mathbb H)$ where $1<p_i< \infty$ for $i = 1, \ldots, m$ then we have
\begin{equation}
\int_{\mathbb R^n} \prod_{j=1}^m |f_j(x,t)|^{p_j} dx \;t^{(m-1)n} \leq C \prod_{j=1}^m
\| f_j \|_{H^{p_j}}^{p_j}.
\end{equation}
\end{prop}

The proof is by induction on $m$, using estimate
\begin{equation}
|f(x)|^p(1-|x|)^n \leq C \| f \|^p_{H^p}, \qquad f \in H^p(\mathbb B),
\end{equation}
and an analogous one for the case of the upper half space.

\begin{lem}\label{shur}
Let $\beta > 0$ and $1 < p <\infty$. Then the operator $S$ defined by
\begin{equation}
Sg(z) = \int_{\mathbb H} \left(\frac{t^n}{|z-\overline w|^{2n}}\right)^{1+\beta} g(w) dm_{\beta n-1}(w), \qquad z \in \mathbb H,
\end{equation}
is bounded on $L^p(\mathbb H, dm_{\beta n-1})$.
\end{lem}

For a similar result for analytic functions in the unit ball see \cite{Zh}.

{\it Proof.} We use Schur's test, see \cite{Zh} for statement and applications to related problems. Let $1/p + 1/q = 1$ and set $h(z) = t^\lambda$, where $\lambda$ will be subject to additional conditions.
Using Lemma \ref{omit} we obtain
\begin{equation}
\int_{\mathbb H} \left(\frac{t^n}{|z-\overline w|^{2n}}\right)^{1+\beta} s^{\lambda q} dm_{\beta n-1}(w) \leq
C t^{\lambda q},
\end{equation}
provided $-\beta n < \lambda q < \beta n + n$. Similarly we have
\begin{equation}
\int_{\mathbb H} \left(\frac{t^n}{|z-\overline w|^{2n}}\right)^{1+\beta} t^{\lambda p} dm_{\beta n-1}(w) \leq
C t^{\lambda p},
\end{equation}
provided $-n - 2\beta n < \lambda p < 0$. Clearly, a parameter $\lambda$ satisfying all the inequalities exists, it suffices
to choose $\max(-\beta n/q, (-n-2\beta n)/p) < \lambda < 0$. Schur's test implies that $S$ is bounded on
$L^p(\mathbb H, dm_{\beta n-1})$ for all $1<p<\infty$. $\Box$

The following result was proved for analytic functions in \cite{LS2} (Theorem 3.3).

\begin{thm}
Let $0<p<q<\infty$, $\alpha > 0$ and let $\mu$ be a positive Borel measure on $\mathbb H$. Then the following two conditions
are equivalent:

$1^o$.
\begin{equation}
\int_{\mathbb H} \left( \int_{\mathbb H} \left( \frac{t^n}{|z-\overline w|^{2n}} \right)^{1+\alpha q} d\mu(z) \right)^{\frac{q}{q-p}} dm_{\alpha qn - 1}(w) < \infty.
\end{equation}

$2^o$.
\begin{equation}
\sum_{k=1}^\infty \eta_k^{-(\alpha qn +n)\frac{p}{q-p}} \mu(\Delta_k)^{\frac{q}{q-p}}<\infty.
\end{equation}
\end{thm}

{\it Proof.} Let us prove implication $1^o \Rightarrow 2^o$. We have
\begin{align*}
\infty & > \int_{\mathbb H} \left( \int_{\mathbb H} \left( \frac{t^n}{|z-\overline w|^{2n}} \right)^{1+\alpha q} d\mu(z) \right)^{\frac{q}{q-p}} dm_{\alpha qn - 1}(w)\\
& \geq C \sum_{k=1}^\infty \int_{\Delta_k} \left( \sum_{j=1}^\infty \int_{\Delta_j} \left( \frac{t^n}{|z-
\overline w|^{2n}} \right)^{1+\alpha q} d\mu(z) \right)^{\frac{q}{q-p}} dm_{\alpha qn - 1}(w)\\
& \geq C \sum_{k=1}^\infty \int_{\Delta_k} \left( \int_{\Delta_k} \left( \frac{t^n}{|z- \overline w|^{2n}} \right)^{1+\alpha q} d\mu(z) \right)^{\frac{q}{q-p}} dm_{\alpha qn - 1}(w)\\
& \geq C \sum_{k=1}^\infty \eta_k^{-(\alpha qn+n)\frac{p}{q-p}} \mu(\Delta_k)^{\frac{q}{q-p}},
\end{align*}
where in the first estimate  we used Lemma \ref{LemmaA} and in last one we used Lemma \ref{LemmaC}.

Next we prove $2^o \Rightarrow 1^o$. Assume $f : \mathbb H \rightarrow \mathbb R$ is subharmonic. Then $|f|^{q/p}$ is
also subharmonic and we have, using Lemma \ref{LemmaB},
\begin{align*}
\int_{\mathbb H} |f(z)| d\mu(z) & = \sum_{k=1}^\infty \int_{\Delta_k} |f(z)|d\mu(z) \leq \sum_{k=1}^\infty
\sup_{z\in\Delta_k} |f(z)| \mu(\Delta_k)\\
& \leq \sum_{k=1}^\infty \left( \int_{\Delta_k^\ast} |f(z)|^{q/p} dm_{\alpha qn -1}(z) \right)^{p/q}
\eta_k^{\frac{-p(n+\alpha qn)}{q}} \mu(\Delta_k).
\end{align*}
Note that $q/p$ and $q/(q-p)$ are a pair of conjugate exponents, hence an application of Holder's inequality gives
\begin{equation}\label{sub}
\int_{\mathbb H} |f(z)| d\mu(z) \leq C \| f \|_{L^{q/p}(\mathbb H,dm_{\alpha qn - 1})} \left[ \sum_{k=1}^\infty \eta_k^{\frac{-(n+qn\alpha)p}{q-p}}
\mu(\Delta_k)^{\frac{q}{q-p}} \right]^{\frac{q-p}{q}}.
\end{equation}
Note that here we used finite overlapping property of the family $\Delta_k^\ast$. Let us choose $g \in L^{q/p}(\mathbb H,
dm_{\alpha qn -1})$, $g \geq 0$, and set
$$f(z) = (Sg)(z) = \int_{\mathbb H} \left( \frac{t^n}{|z-\overline w|^{2n}} \right)^{1+\alpha q} g(w) s^{\alpha qn-1} dw.$$
Since, for any fixed $w \in \mathbb H$, the function $t^n/|z-\overline w|^{2n}$ is subharmonic it easily follows that
$f$ is subharmonic. Moreover, the operator $S$ is bounded on $L^{q/p}(\mathbb H, dm_{\alpha qn-1})$ by Lemma \ref{shur}, i.e. $\| Sg \|_{L^{q/p}(\mathbb H, dm_{\alpha qn-1})} \leq C \| g \|_{L^{q/p}(\mathbb H, dm_{\alpha qn-1})}$.

Now we apply (\ref{sub}) to $f = S(g)$ and obtain

\begin{align*}
& \int_{\mathbb H} \int_{\mathbb H} \left( \frac{t^n}{|z-\overline w|^{2n}} \right)^{1+\alpha q} g(w) s^{\alpha qn-1} dw d\mu(z)\\
\leq & C \| g \|_{L^{q/p}(\mathbb H, dm_{\alpha qn-1})} \left[ \sum_{k=1}^\infty \eta_k^{\frac{-(n+qn\alpha)p}{q-p}}
\mu(\Delta_k)^{\frac{q}{q-p}} \right]^{\frac{q-p}{q}}.
\end{align*}
Since the last estimate is true for any non negative $g \in L^{q/p}(\mathbb H, dm_{\alpha qn -1})$, $1^o$ follows by
duality. $\Box$

An application of the above theorem to $d\mu(z) = |f(z)|^{q-p} dm_\alpha (z)$ gives a characterization of $A(p,q,\alpha)$ spaces for $0<p<q<\infty$, $\alpha > -1$, see \cite{LS2}.

Set $r_k = 2^{-k}$, $I_k = [r_{k+1}, r_k)$ and $H_k = \mathbb R^n \times I_k$ for $k \in \mathbb Z$.
\begin{df}
Let $\mu$ be a positive Borel measure on $\mathbb H$ and let $0<p,q<\infty$. We define the space $K^p_q(\mu)$ as the
space of $f \in h(\mathbb H)$ such that
\begin{equation}
\| f \|_{K^p_q(\mu)}^q = \sum_{k=-\infty}^\infty \left( \int_{H_k} |f(z)|^p d\mu(z) \right)^{q/p} < \infty.
\end{equation}
\end{df}

The spaces $K^p_q(\mu)$ are harmonic Herz type spaces, related classes of spaces appeared in literature, for example in \cite{KY}. Clearly these spaces include weighted Bergman spaces: $K^p_p(m_s) = A(p,s)$, $0<p<\infty$, $s > -1$.

Note that the $K^p_q(\mu)$ (quasi) norm of $f$ is the $l^q$ norm of the sequence $\| f \|_{L^p(H_k, d\mu)}$, $-\infty < k < \infty$. Therefore we have:
\begin{equation}\label{Kpqmon}
\| f \|_{K^p_{q_2}(d\mu)} \leq \| f \|_{K^p_{q_1}(d\mu)}, \qquad 0< q_1 \leq q_2 < \infty.
\end{equation}

\begin{thm}\label{osam}
Let $\alpha > -1$, $0<p\leq q < \infty$ and let $\mu$ be a positive Borel measure on $\mathbb H$. Then the following two conditions are equivalent.

$1^o$. $A(p,\alpha) \hookrightarrow K^p_q(\mu)$.

$2^o$. The measure $\mu$ satisfies the following Carleson type condition:

\begin{equation}
\frac{\mu(\Delta_j)}{|\Delta_j|^{1+ \frac{\alpha}{n+1}}} \leq C, \qquad j \geq 1.
\end{equation}
\end{thm}

{\it Proof.} Sufficiency, for $p=q$, is proved in \cite{AS4}, and the general case follows from the embedding
(\ref{Kpqmon}).

Now we prove necessity, i.e. we assume $1^o$ holds. We choose $l \in \mathbb N_0$ such that $p(n-1+l) > \alpha + n + 1$.
Let us fix a cube $\Delta_j$ and choose $k \in \mathbb Z$ such that $m(H_k \cap \Delta_j) \geq |\Delta_j|/2$. We use a
test function $f = f_{\zeta_j,l}$:
\begin{align}
\| f_{\zeta_j,l} \|_{K^p_q(\mu)}^q & \geq \left( \int_{H_k \cap \Delta_j} |f_{\zeta_j,l}(z)|^p d\mu(z) \right)^{q/p} \geq C \left[ \mu(\Delta_j) \eta_j^{-p(n-1+l)} \right]^{q/p} \notag \\
& = C \mu(\Delta_j)^{q/p} \eta_j^{-q(n-1+l)}. \label{fkpqmu}
\end{align}
Next we estimate $A(p,\alpha)$-norm of $f_{\zeta,l}$ using Lemma \ref{omit}:
\begin{equation}\label{fapa}
\| f \|_{A(p,\alpha)}^p = \int_{\mathbb H} |f_{w,l}(z)|^p t^\alpha dz \leq C \int_{\mathbb H} \frac{t^\alpha dz}{|z-
\overline \zeta_j|^{p(n-1+l)}} \leq C \eta_j^{\alpha + n + 1 - p(n-1+l)}.
\end{equation}
Combining our assumption $\| f\|_{K^p_q(\mu)} \leq C \| f \|_{A(p,\alpha)}$ with (\ref{fkpqmu}) and (\ref{fapa}) gives
$2^o$. $\Box$

\begin{rem}
Sharp embedding theorems of the spaces $K^p_q(\mu)$ into $B^{p,q}_\alpha$ and $K^p_q(\mu)$ into $F^{p,q}_\alpha$
with some restriction on the indexes are also true and they follow from
Theorem \ref{osam} and here our arguments repeat arguments of proof of Theorem
3 which was obtained from linear case (Theorem 2) and appropriate
embedding results  for $B^{p,q}_\alpha$ and $F^{p,q}_\alpha$ spaces.
\end{rem}


\begin{thebibliography}{99}




\bibitem{AS4} M. Arsenovi\'c, R. F. Shamoyan, {\it On embeddings, traces and multipliers in harmonic function spaces},
                  arXiv:1108.5343v1.

\bibitem{AS5} M. Arsenovi\'c, R. F.Shamoyan, {\it Trace theorems in harmonic function spaces on $(\mathbb R^{n+1}_+)^m$ and multipliers theorems}, preprint, 2011, 10 pages.

\bibitem{K1} K. L. Avetisyan, {\it Fractional integro-differentiation in harmonic mixed norm spaces on a half-space},
                  Comment. Math. Univ. Carolinae, Vol. 42 (2001), No. 4, 691-709.

\bibitem{CKY} B. R. Choe, H. Koo and H. Yi, {\it Carleson type conditions and weighted inequalities for harmonic functions}, preprint, 2010.




\bibitem{Co1} W. S. Cohn, {\it Generalized area operators on Hardy spaces}, J. Math. Anal. Appl., 1997, 216: 112-121.

\bibitem{D1} A. E. Djrbashian, {\it The classes $A^p_\alpha$ of harmonic functions in half-spaces and an analogue of
                   M. Riesz' theorem}, Izv. Akad. Nauk Arm. SSR, Matematika 22 (1987), no. 4, 386–398 (in
                   Russian); English transl.: Soviet J. Contemp. Math. Anal. (Armenian Academy of Sciences)
                    22 (1987), no. 4, 74–85.

\bibitem{DS} M. Djrbashian and F. Shamoian, {\it Topics in the theory of $A^p_\alpha$ classes}, Teubner Texte zur
                  Mathematik, 1988, v 105.

\bibitem{GLW} M. Q. Gong, Z. J. Lou, Z. J. Wu, {\it Area operators from $\mathcal H^p$ spaces to $L^q$ spaces}, Sci. China Math., 2010, 53(2): 357-366, doi: 10.1007/s11425-010-0031-9.

\bibitem{Gr1} L. Grafakos, {\it Classical Fourier Analysis}, Graduate Texts in Mathematics 249, Springer 2008.






\bibitem{KY}  H. Koo., K. Nam, H. Yi, {\it Weighted harmonic Bergman kernel on half-spaces}, J. Math. Soc. Japan
             Vol. 58, No. 2, 2006.

\bibitem{LS1} S. Li, R. Shamoyan, {\it On some estimates and Carleson type measure for multifunctional holomorphic spaces in the unit ball}, Bull. Sci. math. 134 (2010) 144-154.

\bibitem{LS2} S. Li, R. Shamoyan, {\it On some properties of analytic spaces connected with Bergman metric ball}, Bull.
              Iran. Math. Soc. Vol. 34 No. 2 (2008), 121-139.

\bibitem{LS3} S. Li, R. Shamoyan, {\it On some extensions of theorems on atomic decompositions of
Bergman and Bloch spaces in the unit ball and related problems}, Complex Variables and Elliptic Equations
Vol. 54, No. 12, December 2009, 1151-1162.

\bibitem{Sh1} N. A. Shirokov, {\it Some generalizations of the Littlewood-Paley theorem}, Zap. Nauch. Sem. LOMI 39 (1974), 162-175; J. Soviet Mat. 8 (1977), 119-129.


\bibitem{OF1}  J. M. Ortega, J. Fabrega, {\it Holomorphic Triebel-Lizorkin Spaces}, Journal of Functional Analysis 151,
                  (1997) 177-212.






\bibitem{St}    E. M. Stein, {\it Singular Integrals and Differentiability Properties of Functions}, Princeton
                   Univ. Press, Princeton, New Jersey, 1970.

\bibitem{StW}  E. M. Stein, G. Weiss, {\it Introduction to Fourier analysis on Euclidean spaces}, Princeton University
                  Press, 1971.

\bibitem{W1} Z. Wu, {\it Area operator on Bergman spaces}, Science in China A (2006), 36(5), 481-507

\bibitem{Zh}   K. Zhu, {\it Spaces of Holomorphic Functions in the Unit Ball}, Springer-Verlag, New York, 2005.

\end{thebibliography}
\end{document}